\documentclass[reqno,final]{amsart}
\usepackage{natbib}  
\usepackage{fancyhdr} 
\usepackage{color} 
\usepackage{hyperref} 
\usepackage{graphicx} 
\usepackage{epstopdf}

\usepackage{amssymb}
\usepackage{graphicx}
\usepackage{bbm}
\usepackage{mathrsfs}
\usepackage{amsmath}
\usepackage{amsthm}
\usepackage{enumerate}
\usepackage{amsfonts,amsopn,amscd}
\usepackage{pstricks}
\usepackage{amssymb}

\definecolor{aleacolor}{rgb}{0.16,0.59,0.78}

\hypersetup{
breaklinks,
colorlinks=true,
linkcolor=aleacolor,
urlcolor=aleacolor,
citecolor=aleacolor}

\renewcommand{\headheight}{11pt}
\renewcommand{\cite}{\citet}

\theoremstyle{plain}
\newtheorem{theorem}{Theorem}[section]                                          
\newtheorem{proposition}[theorem]{Proposition}                          
\newtheorem{lemma}[theorem]{Lemma}

\newtheorem{conjecture}[theorem]{Conjecture}
\theoremstyle{definition}

\theoremstyle{remark}
\newtheorem{remark}[theorem]{Remark}
\newtheorem{example}[theorem]{Example}

\makeatletter \@addtoreset{equation}{section} \makeatother
\renewcommand\theequation{\thesection.\arabic{equation}}
\renewcommand\thefigure{\thesection.\arabic{figure}}
\renewcommand\thetable{\thesection.\arabic{table}}

\newcommand{\aleaIndex}[1]{\href{http://alea.impa.br/english/index_v#1.htm}{\bf #1}}
\eheader{ALEA, Lat. Am. J. Probab. Math. Stat.}{\aleaIndex{12} (1)}{2015}{25}{37}

\elogo{\parbox[c]{81pt}{\includegraphics[width=81pt]{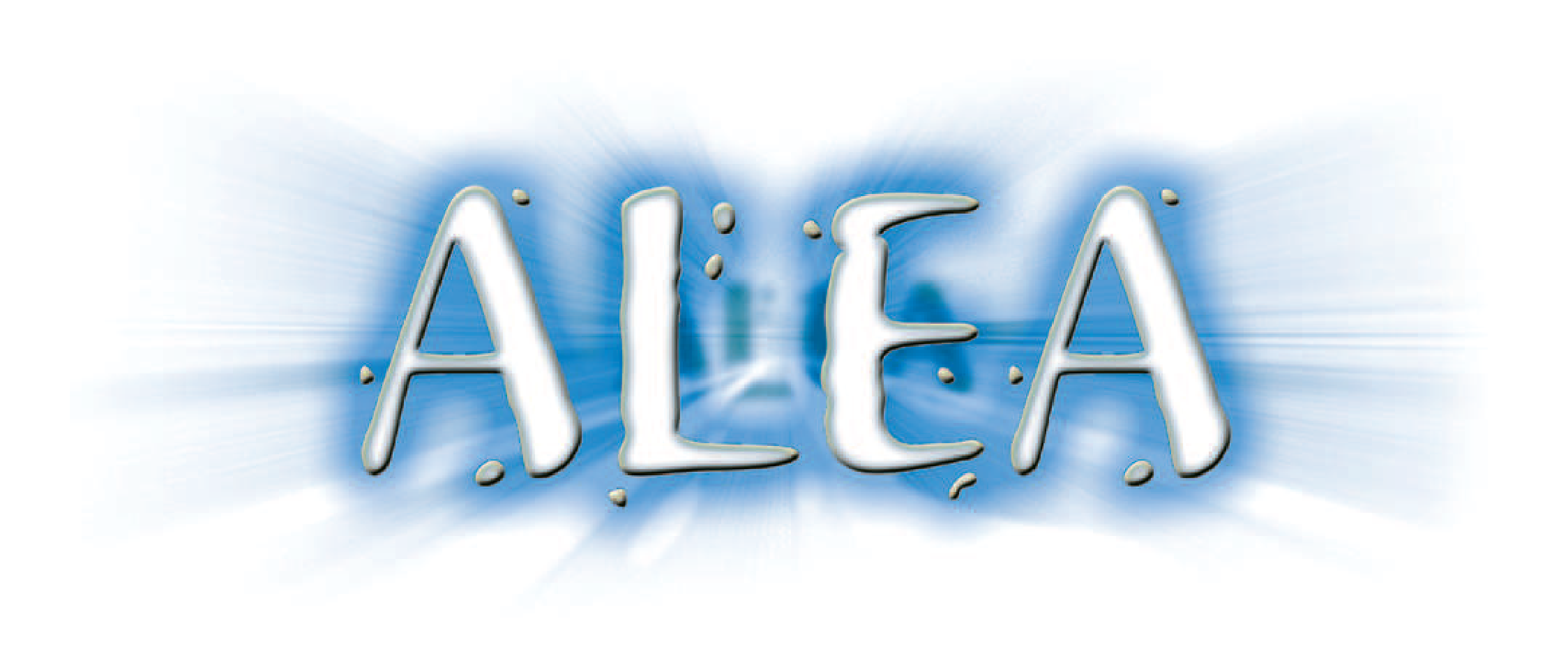}}} 

\renewcommand{\P}{\mathbb{P}}
\newcommand{\E}{\mathbb{E}}
\newcommand{\T}{\mathbb{T}}

\newcommand{\Z}{\mathbb{Z}}
\newcommand{\N}{\mathbb{N}}

\newcommand{\hs}{\hspace{2mm}}
\newcommand{\hsl}{\hspace{1mm}}
\newcommand{\ind}{\mathbbm{1}}
\DeclareMathOperator{\BRW}{BRW}

\date{May 9, 2014; accepted January 19, 2015}

\begin{document}
\title[Number of ends of critical BRW]{The number of ends of critical branching\\ random walks}

\author{Elisabetta Candellero}
\address{University of Warwick\newline
Department of Statistics,\newline
Coventry CV4 7AL, UK.}
\email{elisabetta.candellero@gmail.com}
\urladdr{\url{http://sites.google.com/site/candelleroelisabetta/}}

\author{Matthew I.~Roberts}
\address{University of Bath\newline
Department of Mathematical Sciences,\newline
Bath BA2 7AY, UK.}
\email{mattiroberts@gmail.com}
\urladdr{\url{people.bath.ac.uk/mir20/}}


\subjclass[2010]{60J10, 60J80.} 
\keywords{Branching random walk, trace, topological ends.}

\begin{abstract}
We investigate the number of topological ends of the trace of branching random walk (BRW) on a graph, giving a sufficient condition for the trace to have infinitely many ends. We then describe some interesting examples of non-symmetric BRWs with just one end.
\end{abstract}

\maketitle

\section{Introduction}
Consider a branching random walk (BRW) on a graph described as follows. We are given a locally finite graph $G$ with vertex set $V(G)$ and edge set $E(G)$, and a vertex $v\in V(G)$. We begin with one particle at $v$ at time $0$. For each $n\geq0$, each particle alive at time $n$ gives birth to an independent random number of offspring particles according to some probability distribution $\mu$, each of which independently takes a step according to a specified random walk on the graph with transition kernel $P$. 
Subsequently, the older particle dies.
The resulting configuration of particles comprises the system at time $n+1$.

We write $\BRW(G,\mu,P)$ for such a branching random walk. We assume throughout that $G$ is infinite, $P$ is irreducible, $\mu$ has finite mean $m=m_\mu>1$ and $\mu(0)=0$. We write $(X_n,n\geq0)$ for a random walk on $G$ with transition kernel $P$ ($G$ and $P$ will therefore not be explicit in the random walk notation, but this should not lead to confusion), and $\P_v$ for a probability measure under which the BRW and the random walk are independent and begin at vertex $v$. Let $\rho=\rho(P)=\limsup_{n\to\infty}\P_v(X_n=v)^{1/n}$, the spectral radius of $P$ (which is independent of the choice of $v$).

We say that a BRW is \emph{recurrent} if some vertex of the graph is visited infinitely often by the particles of the BRW with positive probability. Otherwise it is \emph{transient}. Furthermore, a BRW is said to be \emph{strongly recurrent} if all vertices of the graph are visited infinitely many times almost surely (see e.g.\ \citealp[Section 6]{benjamini_peres}).

Benjamini and Peres (see \citealp{benjamini_peres}) showed that if $m<\rho^{-1}$ then the BRW is transient, while if $m > \rho^{-1}$ then the BRW is recurrent. We call a $\BRW(G,\mu,P)$ \emph{critical} if $m_\mu = \rho^{-1}$.

\cite{gantert_mueller} and independently (in a more general setting) \cite{bertacchi_zucca} proved that any critical BRW is transient. The \emph{trace} of the critical process is the subgraph induced by edges that are traversed by particles of the BRW, and is an interesting random structure in its own right. \cite{benjamini2010trace} showed that on Cayley graphs, simple random walk on the trace of transient simple BRW is itself almost surely transient, but any non-trivial simple \emph{branching} random walk on the trace is strongly recurrent.

The range $1<m_\mu \leq \rho^{-1}$ (when non-trivial) leads to the existence of an intermediate phase, in which the process survives globally but eventually vacates every finite set of vertices. Such intermediate behavior has been investigated on different structures, in both discrete and continuous settings. 

In the discrete case, one object of interest is the boundary of the trace of BRW. \cite{hueter_lalley} show that the Hausdorff dimension of the boundary of the trace of a transient BRW defined on a homogeneous tree cannot exceed one half the dimension of the boundary of the tree. This result has been generalized to BRWs defined on free products of groups (see \citealp{candellero_gilch_mueller}).

In a continuous setting, \cite{lalley1997hyperbolic} investigated similar questions for branching Brownian motion on the hyperbolic disk. Their results have been extended to more general Lobachevsky spaces by Karpelevich, Perchersky and Suhov (see \citealp{karpelevich_pechersky_suhov}). Furthermore,  \cite{grigoryan_kelbert} studied the problem of recurrence and transience of branching diffusion processes on Riemannian manifolds. \cite{cammarota_orsingher} studied a branching process where the particles move at finite velocity along the geodesics of the Poincar\'e disk.

An analogous intermediate phase appears also for contact processes on trees (see e.g., \cite{pemantle}, \cite{liggett-cp}, \cite{stacey}, and \cite{liggett-brw-cp}).
For a survey on the behavior of stochastic processes and their characterization of phase transitions on non-amenable graphs, see \cite{lyons-phase-transitions}.

In this work we investigate the number of topological ends of the trace. Roughly speaking, if we remove a large ball about the origin, how many connected components do we see? We will make this definition precise in Section \ref{notation}.

Motivation for investigating this kind of question was outlined in \cite{benjamini_peres}. Classical study of random walks on graphs and groups highlights geometric properties of the space that are captured by the behaviour of the random walk. However, if the graph is large enough then the random walk will typically visit only a tiny part of the structure, and sample path properties of the random walk reflect little of the large-scale geometry. Branching random walks allow us to gain geometric information about the graph on a sample path level.

Recall that $\mu(0)=0$, hence the BRW survives for ever almost surely.
We prove the following theorem, which gives a sufficient condition for the trace to have infinitely many ends. We say that a random walk $X_n, n\geq0$ is {\em quasi-symmetric} if there exists $C$ such that for all vertices $v$ and $w$, and all $n\geq0$,
\[\P_v(X_n=w) \leq C\P_w(X_n=v).\]
In particular, simple random walk on any graph of bounded degree is quasi-{sym\-metric}.

\begin{theorem}\label{vtran}
Suppose that $X_n, n\geq0$ is quasi-symmetric, and that
\[\sum_{n=0}^\infty (n+1)m_\mu^n\P_v(X_n=v) < \infty.\]
Then the trace of $\BRW(G,\mu,P)$ has infinitely many ends almost surely.
\end{theorem}

In particular if $m_\mu < 1/\rho(P)$, then the condition in Theorem \ref{vtran} holds and the trace has infinitely many ends almost surely. On the other hand if $m_\mu > 1/\rho(P)$, then the BRW is recurrent and the trace therefore has as many ends as the underlying graph almost surely (indeed, almost surely the trace \emph{is} the whole graph). We are therefore mainly interested in the critical case $m_\mu = 1/\rho(P)$.

Of course, there are cases which are not covered by Theorem \ref{vtran}. We highlight several examples, finding cases where critical BRW has one end, infinitely many ends, and everything in between.

\section{Preliminaries}\label{notation}

\subsection{The spectral radius and amenability}

Recall that the spectral radius of $P$ is defined to be
\[\rho = \rho(P) = \limsup_{n\to\infty}\P_v(X_n=v)^{1/n}.\]
It is easy to see that this does not depend on the choice of vertex $v$, due to the irreducibility of $P$. See for example \cite{woess2000} for many more details on the spectral radius of a random walk.

In many cases we will consider simple (isotropic) random walk on a graph $G$, which jumps to each of its current neighbours with equal probability, and write $P_G$ for its transition kernel.

\subsection{Ends of a graph}\label{warmup}

A ray is defined to be a semi-infinite simple (i.e.~non-self-intersecting) path within a graph $G$. Say that two rays $R$ and $S$ are equivalent if there exists a third ray that contains infinitely many vertices from each of $R$ and $S$. This defines an equivalence relation, and an \emph{end} of $G$ is an equivalence class of rays. The number of ends of $G$ is then precisely the number of equivalence classes of rays in $G$. For example, $\Z$ has two ends, $\Z^2$ has one end, and $T_3$, the tree in which every vertex has degree $3$, has infinitely many ends.

For $k\geq0$ and $v\in V(G)$, let $B(v,k)$ be the closed ball of radius $k$ (in the graph distance) about $v$. Note that if $e_k$ is the number of connected components of $G\setminus B(v,k)$, then the number of ends of $G$ is equal to $\lim_{k\to\infty}e_k$.

As a warm-up, we prove the following lemma.

\begin{lemma}\label{lemma:G_inf_end}
If $G$ is vertex-transitive and has infinitely many ends, then the trace of any (irreducible, vertex-transitive) BRW  on $G$ will have infinitely many ends almost surely.
\end{lemma}

\begin{proof}
If the BRW is recurrent, then the trace is the whole graph, so the claim is trivially true. We may therefore assume that the BRW is transient.

Fix $k\geq2$ and a vertex $o\in V(G)$.
Recall that for every $n\geq 0$, $B(o,n)$ denotes the closed ball of radius $n$ (in the graph distance) about $o$.
Since $G$ has infinitely many ends, we may choose $n_k$ such that $G\setminus B(o,n_k)$ has at least $k$ pairwise-disjoint components. Let $C_1,\ldots,C_k$ be $k$ of these components. Since $G$ is vertex-transitive, for any $v\in V(G)$ we may take an automorphism $\sigma_v:V(G)\to V(G)$ such that $\sigma_v(o)=v$. Let $C_i(v) = \sigma_v(C_i)$ for each $i=1,\ldots,k$ and $v\in V(G)$.

Suppose that we start with one particle at $o$. There exists $p_k>0$ such that with probability $p_k$, at time $n_k+1$ we have at least one particle in each $C_i$ for $i=1,\ldots,k$, and none of these particles ever has a descendant that returns to $B(o,n_k)$ after time $n_k+1$.

Now, choose $j\geq 2$ and wait until the first time $\tau_j$ that there are at least $j$ particles alive. Call them $u_1,\ldots,u_j$ (if there are more than $j$ then choose the first $j$ according to some arbitrary order) and let $Y_1,\ldots, Y_j$ be their positions. By the strong Markov property, for each $i=1,\ldots,j$, particle $u_i$ has probability $p_k$ of producing a descendant in each of $C_1(Y_i),\ldots, C_k(Y_i)$ at time $\tau_j + n_k+1$, none of whose descendants return to $B(Y_i,n_k)$ after time $n_k+1$. The probability that this does not happen for any of the $j$ particles is therefore $(1-p_k)^j$. As a result,
\[\P(\hbox{this does not happen for any of the $j$ particles for any } j) = 0\]
and therefore
\[\P(\exists v\in V(G) : C_i(v) \hbox{ contains infinitely many points of BRW trace for each } i) = 1.\]
We deduce that the trace has at least $k$ ends almost surely, and since $k$ was arbitrary, infinitely many ends almost surely.
\end{proof}

\subsection{Cartesian products of graphs}\label{cart_sec}

It is natural, at first, to restrict ourselves to vertex-transitive graphs. If $\rho(P_G)=1$ then the critical branching random walk $\BRW(G,1/\rho(P_G),P_G)$ is just a random walk, so we are interested in graphs for which $\rho(P_G)<1$. By Lemma \ref{lemma:G_inf_end}, if a vertex-transitive graph $G$ has infinitely many ends itself, then the trace of any (irreducible, vertex-transitive) BRW will have infinitely many ends too.

Thus we are interested in one-ended, vertex-transitive graphs with $\rho(P_G)<1$. A natural way to construct such graphs is to take the Cartesian product of two infinite, vertex-transitive graphs $G_1$ and $G_2$, at least one of which satisfies $\rho(G_i)<1$. It is a standard result --- see for example \cite{woess2000}, Theorem 4.10 --- that if this holds then $\rho(G_1\times G_2)<1$.

To be precise, for two graphs $G_1$ and $G_2$ (not necessarily vertex transitive), the Cartesian product $G_1\times G_2$ is the graph with vertex set $V(G_1\times G_2) = V(G_1)\times V(G_2)$ and edge set
\[\begin{split}
E(G_1\times G_2) = &\{\{(a,v),(b,w)\} : (a,b)\in E(G_1) \hbox{ and } v=w, \\
& \hbox{ or } a=b \hbox{ and } (v,w)\in E(G_2)\}.
\end{split}\]
Given transition kernels $P_{G_1}$ and $P_{G_2}$ on $G_1$ and $G_2$ respectively, we can then define a random walk on $G_1\times G_2$ by flipping a fair coin at every step: if it comes up heads, then we take a step along an edge inherited from $G_1$ according to $P_{G_1}$, and if it comes up tails, we take a random walk along an edge inherited from $G_2$ according to $P_{G_2}$. We write $P = \frac12 P_{G_1} \oplus \frac12 P_{G_2}$ for the corresponding transition kernel.

For a slightly more thorough treatment of Cartesian products see Section \ref{appendix}.

\section{\texorpdfstring{$T_3\times\Z$}{T\_3 x Z} and lazy random walk}

Our first example is $T_3\times \Z$, where $T_3$ is the tree in which every vertex has degree $3$. 
Recall that $P_{T_3}$ and $P_{\Z}$ denote the transition kernels of simple random walk on $T_3$ and $\Z$ respectively, and consider $P = \frac12 P_{T_3} \oplus \frac12 P_{\Z}$. It is easy to show that $\rho(P_{T_3}) = 2\sqrt2/3$, and of course $\rho(P_{\Z}) = 1$. From this we deduce (see Lemma \ref{CS2}) that
\[\rho(P) = \frac{1}{2}\cdot \frac{2\sqrt2}{3} + \frac{1}{2}\cdot 1 = \frac{\sqrt2}{3} + \frac12 < 1.\]

\begin{proposition}\label{t3z}
If $m_\mu= 1/\rho(P)$, then the trace of $\BRW(T_3\times\Z,\mu,P)$ has infinitely many ends almost surely.
\end{proposition}

\begin{proof}
Consider the projection of the branching random walk onto $T_3$. At each time $n\geq0$, each particle branches independently into a random number of particles with law $\mu$, and each of these particles moves independently according to \emph{lazy} simple random walk on $T_3$; that is, it stays put with probability $1/2$, otherwise it makes a step according to a simple random walk on $T_3$. Call the transition kernel for this walk $L$ (so $L=\frac12 P_{T_3} + \frac12 I$, where $I$ is the identity matrix). Then
\[\rho(L) = \frac{1}{2}\cdot \frac{2\sqrt2}{3} + \frac{1}{2}\cdot 1 = \frac{\sqrt2}{3} + \frac12 = \rho(P),\]
so $\BRW(T_3,\mu,L)$ is a critical branching random walk and therefore transient. We deduce that each copy of $\Z$ is hit only finitely often by particles in $\BRW(T_3\times\Z,\mu,P)$, and since $T_3$ has infinitely many ends, this property is inherited by the trace of the branching random walk by an argument very similar to that 
of Lemma \ref{lemma:G_inf_end}.
\end{proof}

Note that our proof did not require many detailed properties of the two graphs $T_3$ and $\Z$. In fact all that we used was that critical BRW on $T_3$ has infinitely many ends, and that $\rho(P_\Z)=1$.

\section{\texorpdfstring{$T_3\times T_3$}{T\_3 x T\_3}, purple dots, and a proof of Theorem \ref{vtran}}

The natural next question is to consider what happens when our Cartesian product is of two graphs $G_1$, $G_2$ with $\rho(G_1)<1$ and $\rho(G_2)<1$. We note that our previous tactic no longer works for (for example) simple random walk on $T_3\times T_3$. If we look at projections, then we get $\rho(P)<\rho(L)$, and therefore each copy of $T_3$ is hit infinitely often by particles of the critical BRW (see e.g.\ \citealp[Theorem 3.7]{muller2008recurrence} or \citealp[Corollary 3.6]{bertacchi_zucca14}), which tells us nothing about the number of ends of the trace.

We instead look from a different viewpoint: if we start two independent branching random walks from vertices a long way apart, do they meet? Start one BRW from a vertex $v\in T_3\times T_3$, and call this the \emph{red} process; start the other, independent BRW from vertex $w$, and call it the \emph{blue} process. Colour red any vertex that is hit by a particle in the red process, and colour blue any vertex that is hit by a particle in the blue process. Vertices that are coloured both red and blue we call \emph{purple}. We are interested in whether the number of purple vertices is finite or not.

We work in general, aiming for a proof of Theorem \ref{vtran}, but the reader may wish to continue thinking of $T_3\times T_3$. To check that simple random walk on $T_3\times T_3$ satisfies the condition of Theorem \ref{vtran}, we use \cite[Theorem 2(i)]{cartwright_soardi86} together with a result from \cite{cartwright_soardi}, which we translate into our notation in Section \ref{appendix} as Lemma \ref{CS2}. These tell us that if $(X_n,n\geq0)$ is simple random walk on $T_3\times T_3$, then for any vertex $v$, there exists a constant $C_v$ such that
\[\P_v(X_n = v) \sim C_v \frac{\rho^n}{n^3}\]
where $\rho = 2\sqrt2/3$ and we say that $\P_v(X_n=v)\sim f(n)$ if $\P_v(X_n=v)/f(n)\to 1$ as $n\to\infty$ through values such that $\P_v(X_n=v)>0$.

\begin{proposition}\label{purplefinite}
Suppose that $X_n, n\geq0$ is quasi-symmetric, and that
\[\sum_{n=0}^\infty (n+1)m_\mu^n\P_v(X_n=v) < \infty.\]
Take independent red and blue branching random walks started from any vertices $v$ and $w$ respectively. Then the expected number of purple vertices is finite.
\end{proposition}

\noindent
Once we have established this Proposition, we are able to complete the proof of Theorem \ref{vtran}.

\begin{proof}[Proof of Theorem \ref{vtran} (assuming Proposition \ref{purplefinite})]
For every $n\geq 0$ define $N(n)$ to be the set of particles alive in the branching random walk at time $n$. Now fix  $k>0$, and
let $T = \inf\{n\geq0 : |N(n)| \geq k\}$, the first time that we have more than $k$ particles alive (since $m>1$ and $\mu(0)=0$, $T<\infty$ almost surely). Label these particles $1,\ldots, |N(T)|$. For each particle $r$ at time $T$, given its position, its descendants draw out a BRW; call this $\BRW_r$. Note that, by Proposition \ref{purplefinite}, the intersection of the trace of $\BRW_r$ with the trace of $\BRW_s$ is almost surely finite for each $r\neq s$. Thus the particles of $\BRW_r$ almost surely form at least one topological end distinct from the particles of $\bigcup_{s\neq r}\BRW_s$. We deduce that we have at least $k$ ends almost surely. Since $k$ was arbitrary, we must have infinitely many ends almost surely.
\end{proof}

To prove Proposition \ref{purplefinite}, we will use the following well-known result about branching random walks. The proof of this simple version is an easy exercise. Consider a $\BRW(G,\mu,P)$. Let $N(n)$ be the set of particles alive in the branching random walk at time $n$, and for $u\in N(n)$, write $Z(u)$ for its position in $G$.

\begin{lemma}[Many-to-one]\label{manytoone}
For any vertices $v$,$w$ and any $n\geq0$,
\[\E_v[\#\{u\in N(n) : Z(u) = w\}] = m_\mu^n\P_v(X_n = w)\]
where $(X_n, n\geq0)$ is a random walk driven by $P$.
\end{lemma}

We can now prove our key proposition.

\begin{proof}[Proof of Proposition \ref{purplefinite}]
Let $N^R(n)$ be the set of particles in the red BRW at time $n$, and $N^B(n)$ the set of particles in the blue BRW at time $n$. Using independence and then applying Lemma \ref{manytoone},
\begin{align*}
\E & \left [ \# \{ \textnormal{purple vertices}\}\right ] \leq \E\left[\sum_{x\in G} \sum_{k=0}^\infty \sum_{n=0}^\infty \ind_{\{\exists u\in N^R(k), v\in N^B(n) : Z(u)=Z(v)=x\}}\right]\\
&\leq \sum_{x\in G}\sum_{k=0}^\infty\sum_{n=0}^\infty \E[\#\{u\in N^R(k) : Z(u)=x\}]\E[\#\{v\in N^B(n) : Z(v)=x\}]\\
&= \sum_{x\in G}\sum_{k=0}^\infty\sum_{n=0}^\infty m_\mu^k\P_v(X_k = x)\cdot m_\mu^n\P_w(X_n=x).
\end{align*}
But by quasi-symmetry,
\[\sum_{x\in G} \P_v(X_k = x)\P_w(X_n=x) \leq C \sum_{x\in G} \P_v(X_k = x)\P_x(X_n=w) = C\P_v(X_{k+n}=w),\]
so
\[\E\left [ \# \{ \textnormal{purple vertices}\}\right ] \leq C\sum_{k,n}m_\mu^{k+n} \P_v(X_{k+n} = w) = C\sum_n (n+1)m_\mu^n \P_v(X_n=w).\]
Thus whenever this quantity is finite the Proposition, and therefore Theorem \ref{vtran}, holds.
\end{proof}

\vspace{2.5mm}

We remark here that \cite{muller2009dynamical} obtained the same condition for ensuring that BRW on a Cayley graph is \emph{dynamically stable}: that is, if we construct a BRW and then rerandomise each random walk step at rate 1, there is never a time at which the BRW is recurrent.

\vspace{2.5mm}

Relaxing our assumptions even further, we obtain the following even when the underlying random walk is not quasi-symmetric.

\begin{theorem}\label{gen}
Suppose that
\[\sum_{k,n=0}^\infty \rho(P)^{-(k+n)}\P(X_k =X'_n) < \infty\]
where $(X_n, n\geq0)$ and $(X'_n,n\geq0)$ are independent random walks driven by $P$. If $m_\mu = 1/\rho(P) > 1$, then $\BRW(G,\mu,P)$ has infinitely many ends almost surely.
\end{theorem}

\section{One-ended branching random walks, and some open problems}
It is natural to ask whether a converse of Theorem \ref{vtran} might hold. But we already know that a full converse cannot hold: simple random walk on $T_3\times \Z$ satisfies
\[\sum_{n=0}^\infty (n+1)\rho^{-n}\P_v(X_{n}=v) = \infty,\]
but the corresponding critical BRW has infinitely many ends.

Can we, then, construct a non-trivial critical BRW on a graph $G$ with $\rho(P_G)<1$ whose trace is one-ended? If we allow ourselves to bias our random walk in one direction, then the answer is yes. (For simple random walk, the answer is still yes, but the construction is more difficult and we save it for later.) We return to considering $G = T_3\times\Z$, but this time we bias the random walk so that on $\Z$ it is more likely to move in one direction than the other. 
To be precise, let $P_{T_3}$ drive simple random walk on $T_3$, and $P_{\Z}(p)$ drive the random walk on $\Z$ that moves right with probability $p$ and left with probability $1-p$. Let $P(p) = \frac12 P_{T_3} \oplus \frac12 P_{\Z}(p)$.

\begin{proposition}\label{t3zbias}
For $P(p)$ as above with $p\neq 1/2$, the trace of critical BRW on $T_3\times\Z$ has one end almost surely.
\end{proposition}

\begin{proof}
Just as in the proof that the symmetric case has infinitely many ends, this essentially follows from comparing the spectral radius of lazy random walk on $T_3$ with the spectral radius of $P(p)$. Indeed, just as in the proof of Proposition \ref{t3z}, we have $\rho(L) = \sqrt2/3 + 1/2$; but this time $\rho(P(p)) = \sqrt2/3 + \sqrt{p(1-p)} < \rho(L)$. We deduce that each copy of $\Z$ is hit ``often enough'' that there is only one end; the rest of the proof is concerned with making this statement precise. We resort to the use of purple vertices.

Fix a copy of $\Z$; call it $Z_0$. Since $m=1/\rho(P(p)) > 1/\rho(L)$, the lazy BRW on $T_3$ is recurrent and therefore the (non-lazy) BRW on $T_3\times\Z$ hits $Z_0$ infinitely often. Take a realisation of the BRW, and choose two particles. The descendants of the first particle we call \emph{red}, and the descendants of the second particle we call \emph{blue}. Any site in $Z_0$ that is hit by both red and blue particles we colour purple.

We first show that there are infinitely many purple sites almost surely, and then deduce that the trace is one-ended. Let $p_n$ be the probability that a biased RW started from a site in $Z_0$ is in $Z_0$ at time $n$; otherwise defined, $p_n$ is the return probability of the lazy random walk on $T_3$ at time $n$. Since $\rho(L)>\rho(P(p))$, we can choose $k$ such that 
\begin{equation}\label{eq:pk}
p_k > \rho(P(p))^k.
\end{equation}
Fix a red particle in $Z_0$. Call it $u$. Let $S_0^u = \{u\}$, and for $j\geq1$ define
\[S_j^u = \{\text{descendants of particles in $S_{j-1}^u$ that are in $Z_0$ at time $jk$}\}.\]
Let $Y_j^u = |S_j^u|$ for each $j$. Then $(Y_j^u,j\geq0)$ is a Galton-Watson process whose birth distribution $\lambda$ satisfies (by Lemma \ref{manytoone})
\[\E[\lambda] = \rho(P(p))^{-k}p_k > 1.\]
Let $Q^u$ be the event that this Galton-Watson process survives forever, $Q^u = \{Y_j^u \geq1, \hsl \forall j\geq0\}$. Note that the probability that $Q^u$ occurs does not depend on the choice of $u$; let $q = \P(Q^u) > 0$.

We choose a red particle in $Z_0$ via the following algorithm:\\
Let $u_1$ be the first red particle to hit $Z_0$ (in the case of ties, let $u_1$ be any one of the candidates arbitrarily). If $Q^{u_1}$ occurs, then choose $u_1$.
For each $i\geq2$, if $Q^{u_{i-1}}$ does not occur, then let $u_i$ be the first red particle to hit $Z_0$ that is not a descendant of $u_{i-1}$. If $Q^{u_i}$ occurs, then choose $u_i$.

Since each particle's Galton-Watson process is independent of the others, and each has a fixed probability of surviving forever, the algorithm terminates with probability one. Then we have an infinite sequence $v_1,v_2,\ldots$ of red particles in $Z_0$ each of which is at distance at most $k$ from another (our random walk cannot move further than distance $k$ in time $k$). For any particle $v$, let $Y(v)$ be the projection of that particle's position (in $T_3\times\Z$) onto $Z_0$. Since the process is transient, by taking a subsequence and re-ordering if necessary, we may assume that $Y(v_i) > Y(v_{i-1})$ for each $i\geq1$.

Similarly we can construct an increasing sequence of blue particles $w_1,w_2,\ldots$ in $Z_0$ with spacing at most $k$ (where $k$ is such  that (\ref{eq:pk}) is satisfied). From these two sequences of particles, we see that there are infinitely many pairs of red and blue particles within distance $\lceil k/2\rceil$ of each other. Since the RW has a positive probability of stepping upwards $\lceil k/2\rceil$ times in a row, there exists some $\delta>0$ such that each pair has probability at least $\delta$ of generating a purple vertex independently of the others. Thus we must have infinitely many purple vertices.

Now fix any two vertices $v$ and $w$. If we start a red BRW from $v$ and a blue BRW from $w$, then we get infinitely many purple dots almost surely. That means that if $v$ and $w$ are in the trace, then there are almost surely infinitely many distinct paths between $v$ and $w$ within the trace. Thus the probability that there exist two vertices $v$ and $w$ in the trace, that do not have infinitely many distinct paths between them within the trace, is zero. Since the BRW is transient, and so for any $R>0$ with probability one only finitely many particles ever enter $B(o,R)$, we deduce that the trace has only one end.
\end{proof}


On $T_3\times T_3$, the situation is more complicated. We introduce a bias on one of the factors by choosing one direction in the tree to be ``up'' and sending particles in that direction with probability $p$, and the other direction with probability $1-p$. A rigorous description follows.
  
Choose an isometric embedding $\phi:\Z\to T_3$; that is, choose a two-sided infinite path in $T_3$ and label the vertices on that path $\ldots, -2, -1, 0, 1, 2,\ldots$. For each vertex $v$ in $T_3$, we give it a \emph{height} $h(v) = \arg\min_k d(v,\phi(k)) - \min_k d(v,\phi(k))$. If we view $\phi$ as labelling certain vertices in $T_3$, then $h(v)$ is the label of the closest vertex in $\phi(\Z)$ minus the graph distance between $v$ and $\phi(\Z)$. Note that each vertex with height $h$ has two neighbours of heights $h+1$ and one of height $h-1$. Now define a random walk with transition kernel $P_p$ on $T_3$ which jumps from a vertex with height $h$ to its neighbour of height $h-1$ with probability $p$, and to each of its neighbours of height $h+1$ with probability $(1-p)/2$. Let $\tilde P_p = \frac12 P_p \oplus \frac12 P_{T_3}$. Suppose that for each $p$, $\mu_p$ is an offspring distribution such that $m_{\mu_p} = 1/\rho(\tilde P_p)$.

\begin{example}
The trace of $\BRW(T_3\times T_3,\mu_{1/2},\tilde P_{1/2})$ has infinitely many ends almost surely: this follows from exactly the same argument as Proposition \ref{t3z}. The trace of $\BRW(T_3\times T_3,\mu_{1/3},\tilde P_{1/3})$ also has infinitely many ends almost surely, by Theorem \ref{vtran}. On the other hand, if $p>1/2$, then the trace of $\BRW(T_3\times T_3,\mu_p,\tilde P_p)$ has one end almost surely. The proof of this is almost identical to that of Proposition \ref{t3zbias}, except that we have to look at a copy of $\Z$ embedded in $T_3$.

We do not know how many ends the trace of $\BRW(T_3\times T_3, \mu_p, \tilde P_p)$ has when $p\in(0,1/3)\cup(1/3,1/2)$.
\end{example}

\vspace{3mm}

What if we insist that our random walk is simple isotropic random walk? Can we construct a graph on which the corresponding critical BRW is one-ended? Again, the answer is yes.

\begin{example}
Let $\N_0 = \{0,1,2,\ldots\}$, and let $\T$ be the rooted tree in which every vertex has four children, and every vertex except the root has one parent. We say that the root has generation $0$, its children have generation $1$, and so on. Construct a graph $H$ by joining $n\in \N_0$ to every vertex in generation $n$ and every vertex in generation $n+1$ in $\T$, for each $n\geq0$.

Denote by $P_H$ the transition kernel governing simple random walk on $H$. We claim that $\rho(P_H)\in(0,1)$, but if $m_\mu = 1/\rho(P_H)$, then $\BRW(H,\mu,P_H)$ has one end almost surely.

To check that $\rho(P_H)<1$, let $(X_n, n\geq0)$ be simple random walk on $H$. Let $L_k$ be the vertices in generation $k$ of $\T$, together with vertex $k\in \N_0$. Then
\[\P(X_{n+1} \in L_{k+1} | X_n\in L_k) \geq 4/7\]
and
\[\P(X_{n+1} \in L_j | X_n \in L_k) = 0 \hs\hs \hbox{ for } j\not\in \{k-1,k,k+1\}.\]
By coupling with a random walk on $\Z$ that jumps right with probability $4/7$ and left with probability $3/7$, we see that $\rho(P_H)<1$.

On the other hand, starting at the root $o$ of $\T$, by repeatedly jumping to generation $1$ and then back to the root, we get
\[P(X_{2n} = o) \geq (4/5)^n (1/7)^n\]
so $\rho(P_H)>0$.

Showing that the critical BRW is one-ended can be done in a very similar way to the proof of Proposition \ref{t3zbias}. Choose any two particles in the BRW, and colour their descendants red and blue respectively; any site hit by both a red and a blue particle is coloured purple. Clearly a random walk on $H$ eventually hits $\N_0$ almost surely, so $\N_0$ is hit infinitely often by red particles and infinitely often by blue particles. We can ensure (for example by constructing embedded Galton-Watson processes as in the proof of Proposition \ref{t3zbias}) that there are almost surely infinitely many red-blue pairs within distance at most $k$ of each other. Each of these pairs has at least a fixed probability $\delta > 0$ of creating a purple dot (since this time $\N_0$ does not have bounded degree, we have to allow particles to bounce back and forth between $\N_0$ and $\T$). Thus we have infinitely many purple dots almost surely, and since our initial choice of red and blue particles was arbitrary, this guarantees that we have only one end.
\end{example}

\begin{remark}
By gluing together various copies of graphs already constructed, it is easy to construct graphs on which critical BRW has any number of ends. For example, if we glue two copies of $H$ and one copy of $T_3\times T_3$ at a single vertex, then we have a graph on which critical BRW has one end, two ends, or infinitely many ends, each with positive probability.
\end{remark}

\noindent
\textbf{Open problem.} Does there exist a quasi-symmetric random walk on a graph $G$ with $\rho<1$ such that the corresponding critical BRW is one-ended almost surely?

\vspace{2.5mm}

\begin{remark}
An earlier version of this article mentioned a result of \cite{lalley1997hyperbolic} which says that critical branching Brownian motion on the hyperbolic plane has a limit set that does not have full measure. We hypothesised that by analogy one might expect critical branching random walks on co-compact Fuchsian groups (corresponding to tilings of the hyperbolic plane) to have infinitely many ends. Since then \cite{gilch_mueller} have shown that indeed, critical branching random walks on planar hyperbolic Cayley graphs have infinitely many ends.
\end{remark}

Finally, we relay a conjecture from Itai Benjamini.

\begin{conjecture}[Itai Benjamini, private communication]
On any vertex-transitive graph $G$ with $\rho(P_G)<1$, the trace of $\BRW(G,1/\rho(P_G),P_G)$ has infinitely many ends.
\end{conjecture}

\section{Appendix on product random walks}\label{appendix}

Given random walks $X^{(1)},\ldots, X^{(d)}$ on $G_1,\ldots,G_d$ respectively, and real values $\alpha_1,\ldots,\alpha_d\geq0$ with $\alpha_1+\ldots+\alpha_d = 1$, we define the $(\alpha_1,\ldots,\alpha_d)$-\emph{product random walk} $X$ on $G=G_1\times\ldots\times G_d$ by setting
\[\begin{split}
\P_{(v_1,\ldots,v_d)} & (X_1=(w_1,\ldots,w_d)) = \\
& = \begin{cases}\alpha_l\P(X_1^{(l)} = w_l | X_0^{(l)} = v_l) &\hbox{if } w_k = v_k \hs \forall k\neq l, w_l\neq v_l\\
															   \sum_{k=1}^d \alpha_k \P(X_1^{(k)} = v_k | X_0^{(k)} = v_k) &\hbox{if } w_k = v_k \hs \forall k\\
															   0 &\hbox{otherwise}.\end{cases}
\end{split}\]
We write $P=\alpha_1 P^{(1)} \oplus \ldots \oplus \alpha_d P^{(d)}$ for the transition kernel of this random walk. For vertex-transitive graphs, there are at least two natural choices for $\alpha_i$: we could take $\alpha_i = 1/d$ for each $i$, or $\alpha_i = \text{deg}(G_i)/\sum_j\text{deg}(G_j)$ where $\text{deg}(G_j)$ is the degree of an arbitrary vertex in $G_j$. The latter choice corresponds to simple isotropic random walk on $G_1\times\ldots\times G_d$.

The following lemma, due to \cite{cartwright_soardi}, tells us how to translate results about return probabilities on certain graphs into results about return probabilities on their Cartesian products. We recall that the Cayley graph $G(Y)$ of a group $Y$ with generating set $S$ has as its vertex set the elements of $Y$, with two vertices $v,w\in Y$ joined by an edge if $v=ws$ for some $s\in S$. We write that $\P(X_n = v) \sim b_n$ if $\P(X_n=v)/b_n\to 1$ as $n\to\infty$ through values such that $\P(X_n=v)>0$.

\begin{lemma}\label{CS2}
Suppose that $G_1,\ldots ,G_d$ are Cayley graphs with associated random walks $X^{(1)},\ldots, X^{(d)}$. Let $X$ be the $(\alpha_1,\ldots,\alpha_d)$-product random walk on $G=G_1\times\ldots\times G_d$. Fix $v = (v_1,\ldots,v_d)\in V(G)$. Suppose that there exist constants $C_1,\ldots ,C_d>0$ and $a_1,\ldots ,a_d$ such that for each $k=1,\ldots, d$,
\[\P_{v_k}(X^{(k)}_n = v_k) \sim C_j\frac{\rho_k^n}{n^{a_k}}.\]
Then there exists $C$ such that
\[\P_v(X_n = v) \sim C\frac{(\alpha_1\rho_1 +\ldots+\alpha_d\rho_d)^n}{n^{a_1+\ldots+a_d}}.\]
\end{lemma}

Note in particular that the choice of $\alpha_1,\ldots,\alpha_d$ affects the spectral radius in the obvious way, and has no effect on the polynomial terms. Thus the proofs in earlier sections are unaffected by choosing different $\alpha_1,\ldots,\alpha_d$. For example, the results on $T_3 \times \Z$ where we used $\frac12 P_{T_3} \oplus \frac12 P_{\Z}$ hold also for $\frac35 P_{T_3} \oplus \frac25 P_{\Z}$, which corresponds to simple isotropic random walk on $T_3\times \Z$.

\section*{Acknowledgements}
We would like to thank Itai Benjamini for asking us these questions, and for several helpful discussions. We are also grateful to Steven Lalley, Yuval Peres and Yehuda Pinchover for their help. During this project EC received support from the Austrian Academy of Science (DOC-fFORTE fellowship, project number D-1503000014); Microsoft Research; Austrian Science Fund (FWF): S09606, part of the Austrian National Research Network ``Analytic Combinatorics and Probabilistic Number Theory''; and Project DK plus, funded by FWF, Austrian Science Fund (project number E-1503W01230). MR was supported by a CRM-ISM postdoctoral fellowship, McGill University, the University of Warwick, and EPSRC Fellowship EP/K007440/1.

\bibliographystyle{alea2}
\bibliography{2015_02_Candellero}

\end{document}